\documentclass[11pt]{article}
\usepackage{amsmath}
\usepackage{amsfonts}
\usepackage{amsthm}
\usepackage{amssymb}
\usepackage{nicefrac}
\usepackage{graphicx}
\usepackage{enumerate}
\usepackage{subcaption}
\usepackage{fullpage}
\usepackage{float}
\usepackage{color}
\usepackage{algpseudocode}
\usepackage{hyperref}
\usepackage{placeins}

\theoremstyle{plain}
\newtheorem{thm}{Theorem}[section]
\newtheorem{lem}[thm]{Lemma}

\newtheorem{algorithm}[thm]{Algorithm}

\theoremstyle{definition}

\theoremstyle{remark}
\newtheorem{rmk}[thm]{Remark}

\newcommand{\of}[1]{\!\left(#1\right)}

\newcommand{\setR}{\mathbb{R}}
\newcommand{\setN}{\mathbb{N}}

\newcommand{\iu}{\mathfrak{i}}
\DeclareMathOperator{\calO}{{\cal O}}
\DeclareMathOperator{\calF}{{\cal F}}

\newcommand{\diag}[1]{\mathrm{diag}\of{#1}}

\renewcommand{\iu}{i}

\DeclareMathOperator{\curl}{\mathrm{curl}}

{\par\color{blue}\textsubscript{#1}}%
{\par}

{\par\color{red}\textsubscript{#1}}%
{\par}

\newcommand{\domain}{\Omega}

\newcommand{\wei}{\alpha}
\newcommand{\eig}{\omega}
\newcommand{\eigma}{\eig_{\mathrm{max}}}
\newcommand{\eigmi}{\eig_{\mathrm{min}}}
\newcommand{\KryVec}{r}
\newcommand{\Res}{\mathrm{res}}

\DeclareMathOperator{\LH}{\mathrm{span}}

\graphicspath{{./images/}}
\begin{document}
\title{A Krylov Eigenvalue Solver Based on Filtered Time Domain Solutions}
\author{Lothar Nannen\footnote{Institute of Analysis and Scientific Computing, TU Wien, Vienna, Austria, lothar.nannen@tuwien.ac.at} \and 
  Markus Wess\footnote{Institute of Analysis and Scientific Computing, TU Wien, Vienna, Austria, markus.wess@tuwien.ac.at}}

\date{\today}
\maketitle

\begin{abstract}
\noindent This paper introduces a method for computing eigenvalues and eigenvectors of a generalized Hermitian, matrix eigenvalue problem. The work is focused on large scale eigenvalue problems, where the application of a direct inverse is out of reach. Instead, an explicit time-domain integrator for the corresponding wave problem is combined with a proper filtering and a Krylov iteration in order to solve for eigenvalues within a given region of interest. We report results of small scale model problems to confirm the reliability of the method, as well as the computation of acoustic resonances in a three dimensional model of a hunting horn to demonstrate the efficiency.

\end{abstract}

\section{Introduction}
\label{sec:introduction}

Computing eigenvalues and eigenvectors of large scale eigenvalue problems is still a challenging task in applied mathematics. We consider in this paper the computation of eigenpairs $(\eig^2,v)$ to the generalized matrix eigenvalue problem $Sv=\eig^2 M v$ with sparse Hermitian, positive (semi-)definite matrices $S$ and $M$ generated by a finite element discretization of a Laplacian eigenvalue problem. Even in this most simple setting numerical solvers struggle if the matrices become large and if non-extremal or clustered eigenvalues are sought.

We refer to \cite{Saad} for the standard algorithms for such problems. All of the methods therein rely on the basic principle of a simple power iteration, which might be combined e.g. with a Krylov subspace method. Unfortunately, with this basic idea only eigenvalues with largest absolute value can be computed efficiently. One remedy is the use of inverse or shift-and-inverse iterations, which allow in general the computation/approximation of eigenvalues with smallest absolute value or closest to a chosen shift parameter. The price to pay is the need to invert a large, sometimes indefinite matrix in each iteration step. If the matrix dimension is too large for a direct solver to be feasible, then in each step iterative solvers have to be used leading to very high computation times  and/or the need for efficient preconditioners.  

There exist several variants of these methods like the filtered subspace iteration (FEAST), see \cite{FEAST,Gopalakrishnan:20}, or the contour integral method introduced in \cite{Sakurai:03,Beyn} for eigenvalue problems, which are non-linear in the eigenvalue. These methods can reduce the number of iterations, since they focus the iterations to the sought eigenvalues similar to shift-and-inverts described above. 
Again in each iteration step of these methods several linear systems of equations have to be solved. Hence these methods are efficient as long as a direct solver can be used. 

The locally optimal block preconditioned conjugate gradient method (LOBPCG), see \cite{Knyazev:01}, avoids the solution of large linear systems of equations with the additional requirement of a preconditioner for the matrix $S$. With this method, even for large systems, the smallest eigenvalues can be computed efficiently. However, we also consider situations, where not only the smallest eigenvalues but also eigenvalues within a certain range of eigenvalues are of interest. For the latter, the LOBPCG method cannot be used in a straightforward way. 

In this paper we propose a method, which is related to time domain solvers for the Helmholtz problem, see \cite{Grote:20,Appeloeetall:20,Stolk:21}. We construct a Krylov space based on filtered time domain solutions. More precisely, in each Krylov iteration several steps of an explicit time stepping scheme are applied to a semi-discrete wave equation. The time domain solutions are then combined with a proper weight function in order to construct the next Krylov vector. Finally, the large matrix eigenvalue problem is projected and solved on the small dimensional Krylov space. The crucial part of the method is the choice of the weight function. We construct a discrete weight function based on an inverse Fourier transform of a characteristic function of an interval of interest, where the eigenvalues are sought.

The remainder of the paper is organized as follows. The basic concept of the method is presented in Section \ref{sec:basic_example}. Section \ref{sec:numerics} contains numerical examples for a small scale as well as a large scale problem. The small scale problem is used to give hints for choosing suitable method parameters. In the large scale problem we compute resonance frequencies of a three dimensional model of a hunting horn in a closed room. The eigenvalues, which model the playable notes on the horn, belong to a region, where the background eigenvalues of the room are already quite dense, such that the LOBPCG method is not feasible anymore. The same holds true for methods using a direct solver, since the application of such a solver to the systems with more than $10^6$ unknowns typically exceeds the computer memory of a standard desktop computer. The paper closes with a  discussion of extensions to the method in Section \ref{sec:generalizations} and a short conclusion in Section \ref{sec:conclusion}.

\section{Presentation of the method}
\label{sec:basic_example}

In this section we explain the concept in a most simple setting. For a bounded Lipschitz domain $\domain \subset \setR^d$ (for $d=2,3$) we solve for eigenpairs $(\eig^2,u)$ with $\eig\geq 0$ and non-trivial $u \in H^1(\domain)$ of the negative Neumann-Laplacian, i.e., $(\eig^2,u)$ solves
\begin{subequations}\label{def:basicContPDE}
 \begin{align}
   -\Delta u &= \eig^2 u \qquad \text{in } \domain,&\\
   \frac{\partial u}{\partial n} &=0 \qquad \text{at } \partial \domain.&
 \end{align}
\end{subequations}
In the following, we will refer to $\eig$ as resonance or resonance frequency, if $\eig^2$ is an eigenvalue. 

The problem is discretized using a standard Galerkin method. We choose a partition $\mathcal{T}$ of $\domain$ consisting of simplexes and use the discrete finite element space
$$V_h:=\left\{ v\in H^1(\domain): \quad \forall T \in \mathcal{T} \quad v|_T \in \mathcal P_p \right\},$$
where $\mathcal P_p$ denotes the space of polynomials up to degree $p\in\setN$. Thus $V_h$ consists of piecewise polynomials. We solve for discrete eigenpairs $(\eig_h^2,u_h)$ with $\eig_h\geq 0$ and non-trivial $u_h \in V_h$ of the variational formulation
\begin{equation*} 
 \int_\domain \nabla u_h \cdot\nabla \varphi \, dx = \eig_h^2 \int_\domain u_h \, \varphi \, dx\qquad \forall \varphi \in V_h.
\end{equation*}
In the following, we omit the index $h$ since we fix the spatial discretization and focus on the equivalent matrix eigenvalue problem to find eigenvectors $v\in \setR^N\setminus \{0\}$ and eigenvalues $\eig^2\geq 0$ with $N:=\dim V_h \in \setN$ such that
\begin{equation}\label{def:basisMatrixEVP}
 S v = \eig^2 M v.
\end{equation}
The self-adjoint matrices $S=(s_{ij})$, $M=(m_{ij})$ are obtained by choosing basis functions $\psi_n$ such that $ V_h=\LH \{\varphi_1,\dots,\varphi_N\}$  with
$$s_{ij}:=\int_\domain \nabla \varphi_i \cdot \nabla \varphi_j \, dx,\qquad m_{ij}:=\int_\domain \varphi_i \, \varphi_j \, dx,\qquad i,j=1,\dots,N.$$
Clearly, $M$ is positive definite and $S$ positive semi-definite.

\begin{rmk}
  Due to the fact that in the following it will become necessary to apply the inverse of the mass matrix $M$, instead of evaluating the integrals exactly one may approximate them by  numerical integration such that the sparsity pattern of $M$ is more favorable (e.g., $M$ is diagonal if mass lumping is used, see \cite{Cohenetal:01}).
  \label{rmk:mass_lumping}
\end{rmk}

\subsection{Krylov eigenvalue solver}\label{sec:auxProblem}
Let us assume that there exists a matrix $C \in \setR^{N\times N}$ such that the eigenvectors of the matrix eigenvalue problem \eqref{def:basisMatrixEVP} are also eigenvectors of the  auxiliary eigenvalue problem
\begin{equation}\label{def:ArnoldiEWP}
C w= \mu w.
\end{equation}
This includes the possibility, that eigenvectors of \eqref{def:ArnoldiEWP} are linear combinations of eigenvectors to \eqref{def:basisMatrixEVP}. Moreover, let
\begin{equation*}
\mathcal{K}_m(C,\KryVec_0):=\LH\{\KryVec_0,C \KryVec_0,\dots,C^{m-1} \KryVec_0\},
\end{equation*}
be the Krylov space of $C$ with a normalized random starting vector $\KryVec_0\in\setR^N$ and $m\in \setN$. An orthonormal basis $\{b_0,\dots,b_{m-1}\}$ of $\mathcal{K}_m(C,\KryVec_0)$ can be computed iteratively using a Gram-Schmidt orthogonalization, i.e., using $b_0=r_0$ and for each $j=1,\dots,m-1$ we compute $C b_{j-1}$ and  
orthonormalize it with respect to $b_0,\dots,b_{j-1}$. 
Using the projection matrix $B_m=\left(b_0 \dots b_{m-1}\right)\in \setR^{N \times m}$ we project the original eigenvalue problem  \eqref{def:basisMatrixEVP}  onto the Krylov space generated by the matrix $C$: find eigenpairs $\left(\eig_m^2,v_m\right)\in \setR\times \setR^m \setminus\{0\}$ such that
\begin{equation}\label{eq:projevp}
 B_m^\top S B_m v_{m} = \eig_{m}^2 B_m^\top M B_m v_{m}.
\end{equation}
Typically, $m$ is small compared to $N$. Hence, the $m$-dimensional eigenvalue problem \eqref{eq:projevp} can be solved with low computational costs.  

Since \eqref{def:basisMatrixEVP} is a Hermitian eigenvalue problem, following  \cite[Theorem 4.6 and Sec. 6.7]{Saad} we expect convergence of the eigenvectors of \eqref{eq:projevp} towards the eigenspaces to the eigenvalues $\mu$ of \eqref{def:ArnoldiEWP} with largest absolute values. Hence, the projected eigenvalues $\eig_{m}^2$ converge towards those eigenvalues $\eig^2$ of the original problem, for which the corresponding eigenvalues $\mu$ of the auxiliary problem \eqref{def:ArnoldiEWP} have largest absolute value.

Following the considerations above the matrix $C$ has to be chosen in a way that the eigenvalues $\mu$ of \eqref{def:ArnoldiEWP} corresponding to the eigenvalues $\omega^2$ of interest of \eqref{def:basisMatrixEVP} have large absolute values, while the remaining eigenvalues of \eqref{def:ArnoldiEWP} are close to zero. Most standard would be to choose a shift parameter $\rho$ such that $S -\rho M$ is regular and use the shift-and-invert matrix $C:=(S-\rho M)^{-1} M$. It is straightforward to show, that with this choice of $C$ the eigenvectors of \eqref{def:basisMatrixEVP} are identical to the ones of \eqref{def:ArnoldiEWP} and that the correspondence of the eigenvalues is given by $\eig^2=\rho+\mu^{-1}$. Hence, the presented method yields approximations to the squared eigenvalues closest to the shift parameter $\rho$. Unfortunately, using this shift-and-invert technique requires the application of the inverse $(S-\rho M)^{-1}$. In other words, in each Krylov step a linear system of equations for the matrix $S-\rho M$ of dimension $N$ has to be solved. For problems small enough such that a direct solver can be used efficiently, this shift-and-invert approach can be applied with reasonable computational costs.

However, in this paper we are interested in problems where the application of a direct solver is out of reach. Instead of using the shift-and-invert matrix 
in the following Section \ref{sec:filtered_sol} we define the operator $C$ based on filtered time-domain solutions of the underlying wave problem.

  Note that our approach differs from applying the  Arnoldi  method directly to the auxiliary problem \eqref{def:ArnoldiEWP}, since we project the original problem \eqref{def:basisMatrixEVP} onto the Krylov space constructed by the auxiliary problem \eqref{def:ArnoldiEWP}. In our experiments it turned out, that a stopping criterion  for the Krylov iterations is easier to construct for the projected original problem \eqref{eq:projevp}. Moreover, we are interested in $\eig$ and not the auxiliary eigenvalues $\mu$.  For a shift-and-invert method the mapping $\eig^2 \mapsto \mu$ is one-to-one, i.e. $\eig$ can be easily computed if $\mu$ is known. This is not the case for the method based on filtered time-domain solutions.

\subsection{Filtered time-domain solutions}
\label{sec:filtered_sol}
To motivate our specific choice of the auxiliary problem \eqref{def:ArnoldiEWP}, we follow the approach from \cite{Appeloeetall:20}, which was developed for scattering problems. For given $\KryVec\in \setR^N$ let $y(\cdot;\KryVec):[0,\infty) \to \setR^N$ be the solution to the semi-discrete wave problem 
\begin{subequations}\label{eq:waveeq}
  \begin{align}
    M \ddot y (t;\KryVec) &= - S y(t;\KryVec),   \qquad \text{for } t>0,\\
    y(0;\KryVec)&= \KryVec, \qquad    \dot{y}(0;\KryVec)= 0,
  \end{align}
\end{subequations}
where $\dot{y}(\cdot;\KryVec)$ and $\ddot{y}(\cdot;\KryVec)$ denote the first and second time derivative. Since $M,S$ are Hermitian matrices, there exists an orthonormal basis of eigenvectors $v_j \in \setR^N, j=1,\dots,N$ with corresponding eigenvalues  $\eig_{j}^2$  to \eqref{def:basisMatrixEVP} and the unique solution to  \eqref{eq:waveeq} is given by
\begin{equation}\label{eq:solreprswaveeq}
 y(t;\KryVec)=\sum_{j=1}^N \cos(\eig_jt) (v_j^\top \KryVec) v_j.
\end{equation}
For a given piecewise continuous weight function $\wei:[0,\infty) \to \setR$ with compact support  we define the integral operator $\Pi_\wei:\setR^N\to \setR^N$ by
\begin{equation}
  \label{eq:Pidef}
 \Pi_{\wei} \KryVec:=\int_0^\infty \wei(t)y(t;\KryVec) \, dt.
\end{equation}
A discrete version of this integral operator will replace the role of the matrix $C$ from the preceding subsection. The following Lemma quantifies the correspondence of eigenpairs of the initial matrix eigenvalue problem \eqref{def:basisMatrixEVP} and the ones of $\Pi_\wei$.
\begin{lem}\label{lem:aequiEWP}
Let $(\eig^2,v)$ be an eigenpair of \eqref{def:basisMatrixEVP} and the filter function $\beta_\wei:[0,\infty) \to \setR$ be defined by 
\begin{equation}\label{eq:filterfunc}
 \beta_\wei(s):=\int_0^\infty \wei(t) \cos\left(t s\right) dt.
\end{equation}
  Then $(\beta_\wei(\eig),v)$ is an eigenpair of $\Pi_{\wei}$, i.e., $\Pi_{\wei} v=\beta_\wei(\eig) v$. Vice versa, if $(\lambda,v)$ is an eigenpair of $\Pi_{\wei}$, then there exists at least one eigenvalue $\eig^2$ of \eqref{def:basisMatrixEVP} such that $\beta_\wei(\eig)=\lambda$ and $v$ belongs to the sum of eigenspaces of those eigenvalues $\eig^2$ of \eqref{def:basisMatrixEVP} for which $\beta_\wei(\eig)=\lambda$.
\end{lem}
\begin{proof}
 Its straightforward to show, that the solution $y(\cdot;v)$ of \eqref{eq:waveeq} is given by
 $y(t;v)=\cos(\eig t) v$ if $(\eig^2,v)$ is an eigenpair of \eqref{def:basisMatrixEVP}. Hence, the first claim holds by definition of $\beta_\wei$ and $\Pi_\wei$. If $(\lambda,v)$ is an eigenpair of $\Pi_{\wei}$, then the representation \eqref{eq:solreprswaveeq} of the solution to the wave equation yields
 \begin{equation*}
   0=\left(\int_0^\infty \wei(t) \cos\left(\eig_j t\right) dt-\lambda\right) v_j^\top v =\left(\beta_\wei(\omega_j)-\lambda\right)v_j^\top v, \qquad j=1,\dots,N,
 \end{equation*}
  since the eigenvectors $v_j$ form an orthonormal basis of $\setR^N$. For $\eig_j$ with $\lambda\neq \beta_\wei(\eig_j)$ this implies $v_j^\top v=0$. Since $v\neq 0$, there exists at least one $\eig$ with $\lambda= \beta_\wei(\eig)$.
\end{proof}
In other words, the original eigenvalue problem and the eigenvalue problem for $\Pi_\wei$ are somehow equivalent. In particular, if we use $\Pi_\wei$ to construct a Krylov space, this Krylov space approximates sums of eigenspaces to the original problem only. If this were not the case Krylov steps would possibly be wasted into approximations of eigenvectors to $\Pi_\wei$ which are irrelevant for the original eigenvalue problem.

In what follows, we motivate our choice of the weight function $\wei$.
Using the symmetric extension $\hat \wei(t):=\wei(-t)$ for $t<0$, the function $\beta_\wei$ can be represented by the Fourier transform of $\hat \wei$:
 \begin{equation*}
  \beta_\wei(s)=\int_0^\infty \wei(t) \cos\left(t s\right) dt =\frac{1}{2} \int_{-\infty}^\infty \hat \wei (t) \exp(-\iu t s) \, dt = \sqrt{\frac{\pi}{2}} \calF(\hat \wei)  \left(s\right).
 \end{equation*}
Since the presented Krylov method converges towards the eigenvalues with largest absolute values, the weight function $\wei$ should be adapted to the location of the sought eigenvalues. E.g., if eigenvalues $\eig$ in an interval $[\eigmi,\eigma]$ are sought, it would be optimal to find $\wei$ such that $\beta_\wei=\chi_{[\eigmi,\eigma]}$ where $\chi_A$ denotes the characteristic function of a set $A$. This would correspond to choosing $\wei$ as 
\begin{equation}
  \label{eq:filterf}
  \sqrt{\frac{2}{\pi}} \calF^{-1}\left( \chi_{[\eigmi,\eigma]} \right)=\frac{4}{\pi t}\sin\left(\frac{t}{2}\left(\eigma-\eigmi\right)\right) \cos \left(\frac{t}{2}\left(\eigma+\eigmi\right)\right),
\end{equation}
for $t>0$.
Obviously, the function on the right hand side is not compactly supported in $[0,\infty)$. However we are not forced to use exactly this weight function. In fact, in the following we choose a finite time interval $[0,T]$ with $T>0$ together with the weight function 

\begin{align}
  \wei(t) := \begin{cases}
    \frac{2(\eigma-\eigmi)}{\pi},&t=0,\\
    \frac{4}{\pi t}\sin\left(\frac{t}{2}\left(\eigma-\eigmi\right)\right) \cos \left(\frac{t}{2}\left(\eigma+\eigmi\right)\right),&t\in(0,T],\\
    0,&t>T.
  \end{cases}
  \label{eq:weightf}
\end{align}
Note that the method is not limited to this specific weight function (see also Section \ref{sec:gen_weight}). In general, for fast convergence the function $\wei$ should be chosen such that $ |\beta_\wei(\eig)| \gg |\beta_\wei(\tilde \eig)|$ for the sought $\eig$ and the unsought $\tilde \eig$.

\subsection{Discretization of the filtered time-domain solution}
We discretize the integral operator \eqref{eq:Pidef} using a rectangle rule. To this end, we introduce a time-stepping method for the wave equation \eqref{eq:waveeq}. Since we are interested in eigenvalue problems where applying a direct solve is out of reach we focus on explicit methods. For a fixed, uniform step-size $\tau>0$ we approximate $y(\ell \tau;\KryVec)$ for $\ell \in \setN$ by $y_\ell(\KryVec)$ using finite differences. This leads to the explicit, second order two-step method 
\begin{equation}\label{eq:2stepmethod}
  y_{\ell+1}(\KryVec)=- \tau^2 M^{-1}Sy_{\ell}(\KryVec)+2y_{\ell}(\KryVec)-y_{\ell-1}(\KryVec),\qquad \ell \in \setN,
\end{equation}
with initial time steps 
\begin{align}
  y_{-1}(\KryVec)= y_{0}(\KryVec)=\KryVec.
\end{align}
Note, that we define $y_{-1}(\KryVec)=\KryVec$ since $y_{-1}(\KryVec)$ approximates $y(-\tau;\KryVec)=y(0;\KryVec)-\tau \dot y (0;\KryVec)+\calO(\tau^2)= \KryVec+\calO(\tau^2)$.

Using the rectangle rule with $L\in \setN$ quadrature points and therefore the step-size $\tau:=\nicefrac{T}{L}$, we finally arrive at the fully discrete linear mapping $C:\setR^N \to \setR^N$:
\begin{equation}\label{eq:disc_filtercoeff}\KryVec_{j-1}\mapsto \KryVec_j=C \KryVec_{j-1}:=\sum_{\ell=0}^{L-1}\tau \wei(\ell\tau) y_\ell(\KryVec_{j-1}),
\end{equation}
which approximates the integral operator \eqref{eq:Pidef}.

In order to find the discrete analogon to the filter function $\beta_\wei$ defined in \eqref{eq:filterfunc} we study the scalar linear difference problems: for $\eig>0$ let $q_{\ell}(\eig)\in\setR$ for $\ell\in\setN_0$ be the unique solution to
  \begin{align}\label{def:scalartimestep}
    q_{-1}(\eig)&=1,&q_0(\eig)&=1,&
    q_{\ell+1}(\eig)&=\left(2-\tau^2\eig^2\right) q_\ell(\eig)-q_{\ell-1}(\eig).
  \end{align} 
Note, that diagonalization of the Matrix $M^{-1} S$ separates \eqref{eq:2stepmethod} into such problems for those $\eig>0$ such that $\eig^2$ are eigenvalues of the matrix eigenvalue problem \eqref{def:basisMatrixEVP}.
It is straightforward to show that
  $$q_{\ell}(\eig)=\frac{1-\overline{c(\eig)}}{2\iu \Im(c(\eig))}c(\eig)^\ell+\frac{c(\eig)-1}{2\iu \Im(c(\eig))}\overline{c(\eig)}^\ell=\frac{\Im(c(\eig)^\ell)-|c(\eig)|^2\Im(c(\eig))^{\ell-1}}{\Im(c(\eig))},$$
  with $$c(\eig):=1-\frac{\tau^2 \eig^2}{2}+ \iu \sqrt{1-\left(1-\frac{\tau^2\eig^2}{2}\right)^2.}$$
  Hence, $|c(\eig)|=1$ if the CFL condition $\tau<\nicefrac{2}{\eig}$ is satisfied. In this case $\left|q_{\ell}(\eig)\right|$ is bounded in $\ell$. Otherwise the sequence $q_{\ell}(\eig)$ is unbounded for $\ell\to\infty$.
  
 \begin{lem}\label{lem:fullydiscfilter}
  The fully discrete filter function $\tilde\beta_\wei:[0,\infty) \to \setR$ defined by
  \begin{equation}\label{eq:fullydiscfilter}
     \eig\mapsto \tilde\beta_\wei(\eig):= \sum_{\ell=0}^{L-1} \tau \wei(\tau \ell)q_{\ell}(\eig).
  \end{equation}
  is a polynomial of degree $L-1$ in $\eig^2$ and there holds
  \begin{equation}\label{eq:polyFilter}
 C=\tilde\beta_\wei\left(M^{-1} S \right).
  \end{equation}
  In particular, if  $(\eig^2,v)$ is an eigenpair of \eqref{def:basisMatrixEVP}, then $(\tilde\beta_\wei(\eig),v)$ is an eigenpair of $C$, and if $(\lambda,v)$ is an eigenpair of $C$ then there exists at least one eigenvalue $\eig^2$ of \eqref{def:basisMatrixEVP} with $\lambda=\tilde\beta_\wei(\eig)$ and $v$ belongs to the sum of eigenspaces of those eigenvalues $\eig^2$ of \eqref{def:basisMatrixEVP} for which $\tilde\beta_\wei(\eig)=\lambda.$
 \end{lem}
 \begin{proof}
   By definition of the scalar time-stepping \eqref{def:scalartimestep} $q_\ell(\eig)$ is a polynomial of degree $\ell$ in $\eig^2$. Hence, $\tilde\beta_\wei(\eig)$ defined in \eqref{eq:fullydiscfilter} is a polynomial of degree $L-1$ in $\eig^2$. To prove \eqref{eq:polyFilter} we start with an arbitrary vector $\KryVec$ and the orthonormal basis of eigenvectors $v_j$ to \eqref{def:basisMatrixEVP} corresponding to the eigenvalues $\eig_j^2$.  The solutions $y_\ell(\KryVec)$ to \eqref{eq:2stepmethod} are given by
  \begin{equation*}
   y_\ell\left(\KryVec\right)=y_\ell\left(\sum_{j=1}^N \left(\KryVec^\top v_j \right) v_j \right)=\sum_{j=1}^N \left(\KryVec^\top v_j \right) y_\ell \left(v_j\right)=\sum_{j=1}^N \left(\KryVec^\top v_j \right) q_\ell(\eig_j) v_j.
  \end{equation*}
  Hence, by definition of the matrix application $C$ in \eqref{eq:disc_filtercoeff} there holds
  \begin{equation*}
  \begin{aligned}              
   C \KryVec&=\sum_{\ell=0}^{L-1} \tau \wei(\ell\tau) \sum_{j=1}^N \left(\KryVec^\top v_j \right) q_\ell(\eig_j) v_j =  \sum_{j=1}^N \tilde\beta_\wei(\eig_j) \left(\KryVec^\top v_j \right) v_j\\
   &=V \diag{\tilde\beta_\wei(\eig_1),\dots,\tilde\beta_\wei(\eig_N)} V^\top \KryVec=\tilde\beta_\wei\left(M^{-1} S\right)\KryVec,
   \end{aligned}
  \end{equation*}
   where we used the spectral decomposition $M^{-1} S=V^\top \diag{\eig_1,\dots,\eig_N} V$ with $V:=\left(v_1,\dots,v_N\right)$. The equivalence of eigenpairs of \eqref{def:basisMatrixEVP} to eigenpairs of $C$ is a direct consequence of \eqref{eq:polyFilter}.
 \end{proof}
In other words, we use a specific type of polynomial filtering. Note, that we are only interested in values of the filter function $\tilde\beta_\wei$ for which the time-stepping with stepsize $\tau$ is stable, since no larger discrete eigenvalues exist.

\begin{rmk}\label{rmk:timestep}
  Replacing the explicit time-stepping \eqref{eq:2stepmethod} in the results above by other types of time-steppings is straightforward as long as the time-stepping is diagonalized by the eigenvectors $v_j$, i.e., if  $y_\ell(v_j)=f_\ell(\eig_j) v_j$ holds for some functions $f_\ell$.
\end{rmk}

\begin{figure} 
\centering
  \includegraphics{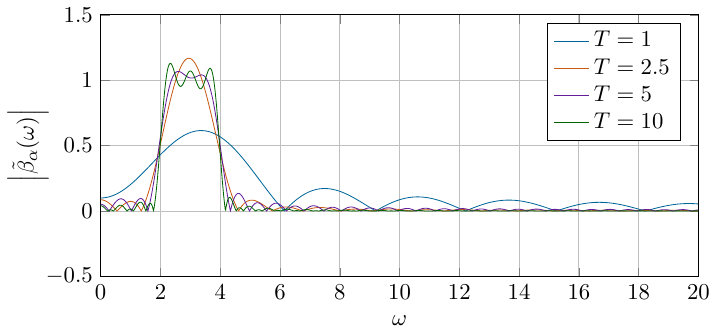}
  \caption{Discrete filter functions for fixed time-step $\tau=0.025$, the target interval $[\eigmi,\eigma]=[2,4]$, and varying end times $T$.}
  \label{fig:filter_4_16}
\end{figure}

Figure \ref{fig:filter_4_16} shows discrete filter functions $\tilde\beta_\wei$ of \eqref{eq:fullydiscfilter} for fixed time-step $\tau$ and $[\eigmi,\eigma]=[2,4]$ in \eqref{eq:filterf}, while varying the end time $T$ (and thus also the number of total time-steps).
While the filter certainly separates the wanted from the unwanted eigenvalues better for larger end times $T$ the total number of necessary time-steps in the Arnoldi process can not be predicted by solely looking at the discrete filter function (see numerical experiments in Section \ref{sec:numerics}).

\subsection{Algorithm}  
The considerations above lead to the following basic algorithm for computing the Krylov subspace $\mathcal{K}_m(C,\KryVec_0)$.
 \begin{algorithm}[Krylov subspace by filtered time-domain solutions] ~\\
   \textbf{Input:} matrices $M^{-1},S\in\setR^{N\times N}$, random starting vector $\KryVec_0\in\setR^N$ with $\|\KryVec_0\|_2=1$, time step size $\tau>0$, number of time steps $L\in \setN$, dimension of Krylov space $m\in\setN$, weight function $\wei$
   \begin{algorithmic}[1]
   \For{ $k=0,\ldots, m-1$} \qquad \hfill Krylov loop
     \State{$y_0 := y_1 := \KryVec_k$} \qquad\hfill initial values for time loop
     \State{$\KryVec_{k+1} := \tau \wei(0) y_0$} \qquad\hfill first term in \eqref{eq:disc_filtercoeff}
     \For{ $i=1,\ldots, L-1$} \qquad \hfill time loop
     \State{$y_2:=-\tau^2M^{-1} S y_1+2y_1-y_0$} \qquad \hfill time iteration \eqref{eq:2stepmethod}
     \State{$\KryVec_{k+1}+=\tau \wei(i\tau) y_2$} \qquad \hfill additions in \eqref{eq:disc_filtercoeff}
     \State{$y_0 := y_1$, $y_1:=y_2$} \qquad \hfill preparation for the next time step     
     \EndFor
   \For{$j=0,\ldots,k$} \qquad \hfill orthogonalization loop
     \State{$\KryVec_{k+1}-=\left(\KryVec_{k+1},\KryVec_j\right)_2 \KryVec_j$}     
    \EndFor
    \If{$\left\|\KryVec_{k+1}\right\|_2\neq 0$}
    \State{$\KryVec_{k+1}*=1/\left\|\KryVec_{k+1}\right\|_2$} \qquad \hfill normalization
    \Else\qquad \hfill exact eigenspace found
    \State{stop Krylov loop} 
    \EndIf
    \EndFor
 \end{algorithmic}
   \textbf{Output:} projection matrix $B_m:=\left(\KryVec_0,\dots,\KryVec_{m-1}\right)\in \setR^{N\times m}$ with orthonormal column vectors.
 \end{algorithm}
The projected eigenvalue problem \eqref{eq:projevp} with small dimension $m$ can be solved using a standard eigenvalue solver. Note, that the main costs of the method are related to the total number of time-steps, i.e., the dimension $m$ of the Krylov space times the number of time steps $L$ in each Krylov iteration. The orthonormalization of Krylov vectors and the solution of the small eigenvalue problems can be neglected for $m\ll N$.

It remains to discuss a stopping criterion, i.e., how to choose the dimension $m$ of the Krylov space, and a criterion to distinguish converged from non-converged eigenvalues. To this end, we compute for the eigenpairs $\left(\eig^2_{j;m},v_{j;m}\right)$ of \eqref{eq:projevp} the residuals
\begin{equation}\label{eq:res}
 \Res_j:=\left\|\left(S-\eig_{j;m}^2 M \right) B_m v_{j;m}\right\|_2,\qquad j=1,\dots,m,
\end{equation}
of the large eigenvalue problem \eqref{def:basisMatrixEVP}. We accept eigenpairs $\left(\eig^2_{j;m},B_m v_{j;m}\right) \in \setR\times \setR^N$ for which the residuals are below a given tolerance. A possible stopping criterion would be to increase $m$ until a fixed number of accepted eigenvalues is found.

Other approaches turned out to be problematic in practice. One idea would be to use the absolute values of the eigenvalues $\mu_{j,m}:=\tilde \beta_\wei(\eig_{j;m})$ of the auxiliary problem \eqref{def:ArnoldiEWP}, which was used to construct the Krylov space. Since we expect convergence towards the eigenvalues $\mu$ with largest absolute values, we e.g., could accept those eigenpairs $\left(\eig^2_{j;m},B_m v_{j;m}\right)$ for which the ratio $\left(\nicefrac{\left|\hat \mu\right|}{\left|\mu_{j,m}\right|}\right)^m$ with $\left|\hat \mu\right|:=\min\left\{\left|\mu_{j,m}\right|, j=1,\dots,m\right\}$ is small enough. 

Nevertheless, in our numerical experiments sometimes non-converged eigenvalues with large $|\mu|$ appeared. This might be the case, if eigenvalues within a certain range of frequencies are sought. Then, from time to time, a discrete eigenvalue, which will finally converge for $m\to N$ from above to an eigenvalue below the region of interest, passes by and leads by accident to a large value $|\mu|$. See Fig.\ref{fig:horn_mid_fq} in the numerics section for such an example, where even after 50 Krylov steps there are non-converged eigenvalues with large values $\mu$. Therefore, we refrained from using the values $\tilde  \beta_\wei(\eig_{j;m})$ in an error indicator. 

Another idea would be to use a standard Arnoldi solver for the auxiliary problem \eqref{def:ArnoldiEWP}, i.e., project the matrix $C$ to the Krylov space, solve the projected eigenvalue problem, use e.g., the residuals of this auxiliary problem as an error indicator, and project the original eigenvalue problem \eqref{def:basisMatrixEVP} only on the span of those eigenvectors, which were computed and accepted by the auxiliary problem. We refrained from this approach since the residuals of the auxiliary problem are somewhat artificial. The eigenvectors are approximations of the sought eigenvectors, but the eigenvalues $\mu$ have no one-to-one relation to the eigenvalues $\eig^2$ or the original eigenvalue problem.

\FloatBarrier

\section{Numerical experiments}
\label{sec:numerics}


We apply our algorithm to two different sets of problems. In Section \ref{sec:small_scale_numerics} we tackle a two-dimensional problem where reference solutions can be computed using a shift-and-invert Arnoldi method. In Section \ref{sec:large_scale_numerics} we choose a three-dimensional problem with more than $10^6$ unknowns, where the application of reference solutions would require significant computational effort when a direct inverse is used in the process. All numerical examples where carried out using the high-order finite element library Netgen/NGSolve (\cite{netgen,ngsolve:14}).

\subsection{Small scale (2d) examples}
\label{sec:small_scale_numerics}
The examples in this subsection are chosen in a way that the resulting systems are small enough to allow to compute reference values by applying a shift-and-invert Arnoldi algorithm. The goal of this subsection is to demonstrate the applicability and functionality of our method as well as to give hints on how to choose the parameters.
\subsubsection{Description of the experiments}
\label{sec:dumbbell_desc}
\begin{figure}
  \centering
  \includegraphics{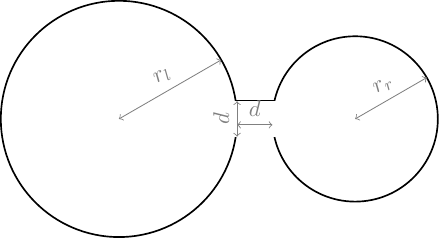}
  \caption{Sketch of the two-dimensional geometry used for the small scale examples from Section \ref{sec:small_scale_numerics}}
  \label{fig:dumbbell_geos}
\end{figure}

\begin{figure}
  \centering
  \includegraphics{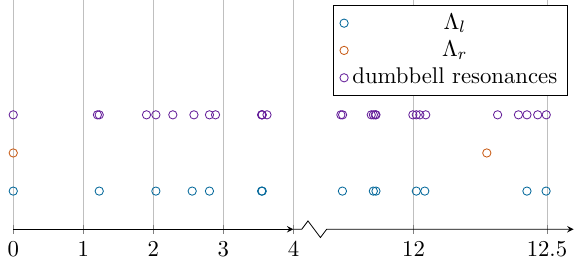}
  \caption{Exact resonances of \eqref{def:basisMatrixEVP} on the circles with radii $r_l=1.5,r_r=0.15$ and the whole dumbbell domain (cf. Figure \ref{fig:dumbbell_geos} and \eqref{eq:dumbbell_params}).}
  \label{fig:dumbbell_evs}
\end{figure}

We choose a two-dimensional domain shaped like a dumbbell, consisting of two circles with radii $r_l,r_r$ connected by a square with width $d$ (cf. Figure \ref{fig:dumbbell_geos}). We prescribe homogeneous Neumann boundary conditions and 
 choose the geometry parameters 
\begin{align}
  r_l &= 1.5,&r_r &= 0.15,&d&=0.03.
  \label{eq:dumbbell_params}
\end{align}
For the spacial discretization we employ the mass lumping technique described in Remark \ref{rmk:mass_lumping}. We choose a mesh-size of $h=d=0.03$ (cf. Figure \ref{fig:2d_ef_low}) and second order finite elements.
We expect the corresponding discrete resonances $\eig_j$ to be perturbations of the union of the three sets
\begin{align*}
  \Lambda_l&=\left\{\frac{\lambda_{j,n}}{r_l},j,n\in\setN_0\right\},&\Lambda_r&=\left\{\frac{\lambda_{j,n}}{r_r},j,n\in\setN_0\right\},&\Lambda_d &= \left\{\frac{\pi \sqrt{n^2+j^2}}{d},j,n\in\setN_0\right\},
\end{align*}
where $\lambda_{j,n}$ are the roots of the derivatives $J_n'$ 
of the cylindrical Bessel functions $J_n$ of the first kind with order $n\in\setN_0$ (cf. \cite[(10.2.2)]{NIST:DLMF}).
Figure \ref{fig:dumbbell_evs} shows the resonances $\eig$ of the discrete problem \eqref{def:basisMatrixEVP} together with the sets $\Lambda_l$ and $\Lambda_r$. The set $\Lambda_d\setminus\{0\}$ is omitted since it consists of high frequencies, which are irrelevant for our computations. In our following experiments we look for the discrete resonances which correspond to perturbations of the smallest non-trivial elements of the sets $\Lambda_l,\Lambda_r$.

\subsubsection{First results}
\label{sec:small_scale_first}

\begin{figure}[t]  
  \centering
  \includegraphics[width=\textwidth]{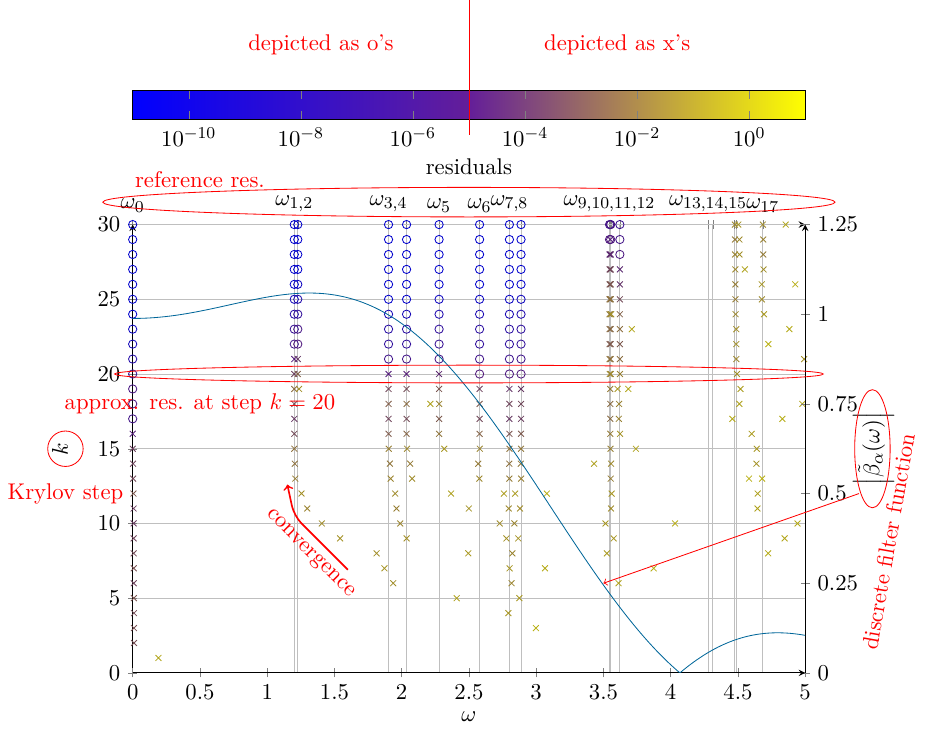}
    \caption{Numerical results for the computation of resonances of \eqref{def:basisMatrixEVP} for the dumbbell domain and discrete filter function for $T=300\tau$, $\eigmi = 0,\eigma=3$.}
  \label{fig:dumbbell_td_evs_qual}  
  \end{figure}
  
\begin{figure}[t]      
  \centering
    \includegraphics[width=0.8\textwidth]{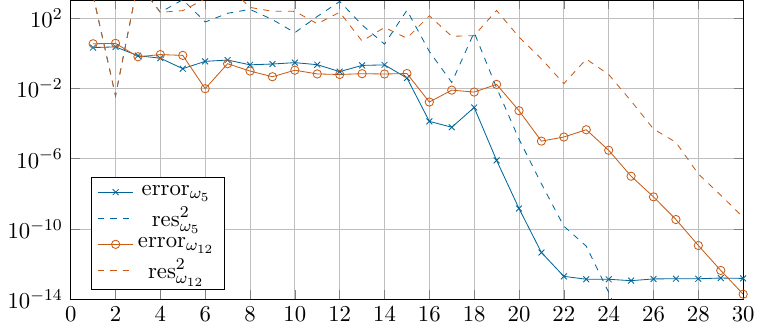}
    \caption{Errors and residuals for selected resonances of Fig. \ref{fig:dumbbell_td_evs_qual}.}
  \label{fig:dumbbell_td_evs_conv}
\end{figure}

A stable time step $\tau\approx 0.0056$ is determined using a power iteration and the CFL condition of the time-stepping \eqref{eq:2stepmethod}. We choose the weight function $\alpha$ given in \eqref{eq:weightf} with $\eigmi=0$, $\eigma=3$ and end-time $T=300\tau\approx 1.68$. Figure \ref{fig:dumbbell_td_evs_qual} shows the resulting approximations to the discrete resonances in each Krylov step $k$, as well as the residuals (see \eqref{eq:res}), reference values and the discrete filter function $\tilde\beta_\alpha$.
Since we will use this type of illustration throughout the remainder of the paper we use Figure \ref{fig:dumbbell_td_evs_qual} to explain all the features of this illustration in more detail.

The horizontal axis corresponds to the values of the approximated resonances $\omega$. The vertical axis on the one hand, corresponds to the Krylov step $k$ and on the other hand to the absolute value of the discrete filter function $\tilde\beta_\alpha(\omega)$. The reference resonances are marked by the vertical grid, thus we observe convergence of the computed values towards these horizontal lines going upward (i.e., increasing the Krylov step). The colors denote the value of the residuals (color bar above the plot), where we use the ring symbol for residuals smaller than $10^{-5}$ (i.e., we consider the resonance converged) and the cross symbol for residuals larger than $10^{-5}$. Figure \ref{fig:dumbbell_td_evs_conv} confirms that our algorithm converges for two selected resonances and that the squares of the residuals behave asymptotically as the errors.

As to be expected, in general, resonances where the value of the filter function is larger are approximated after fewer Krylov steps: Resonances in the range of $0$ to $3$ are approximated well at around $20$ Krylov steps. The resonances around $3.5$ are approximated later. The resonances at around $4.5$ still seem to be approximated (but are by far not converged yet). The resonances close to the root of the filter function around $4.3$ are not approximated at all using this filter and $30$ Krylov steps.

We also note that the convergence does not correspond exactly to the absolute value of the filter function. This is a result of the fact that we expect the convergence of resonances with close absolute values of $\tilde\beta_\wei$ to be worse than the convergence of resonances with more isolated absolute  values of $\tilde\beta_\wei$.

\subsubsection{Different filters}
\label{sec:small_scale_filters}
\begin{figure}
  \begin{subfigure}{0.5\textwidth}
    \includegraphics[width=\textwidth]{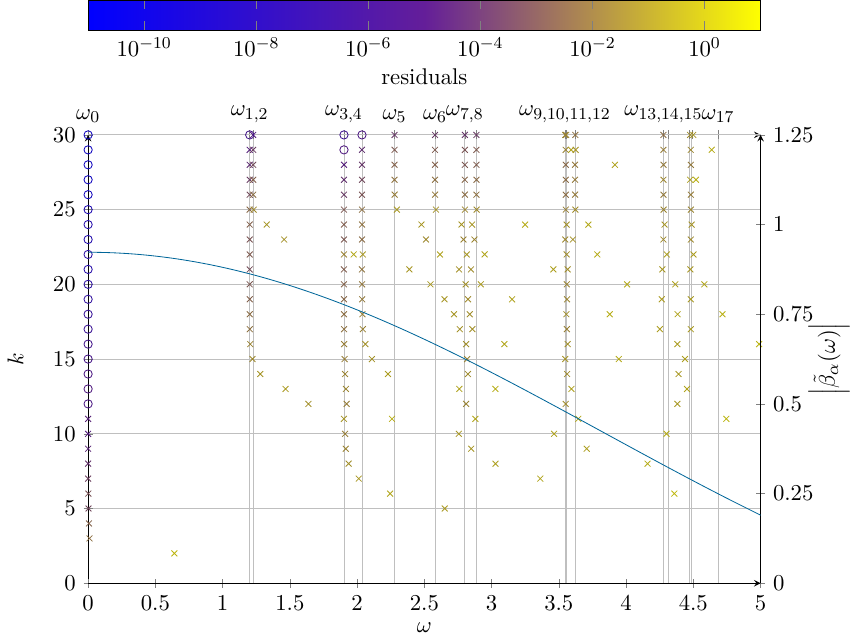}
  \caption{$\eigmi=0,\eigma=3, T=100\tau$}
  \label{fig:dumbbell_0_3_100}
  \end{subfigure}
  \begin{subfigure}{0.5\textwidth}
  \includegraphics[width=\textwidth]{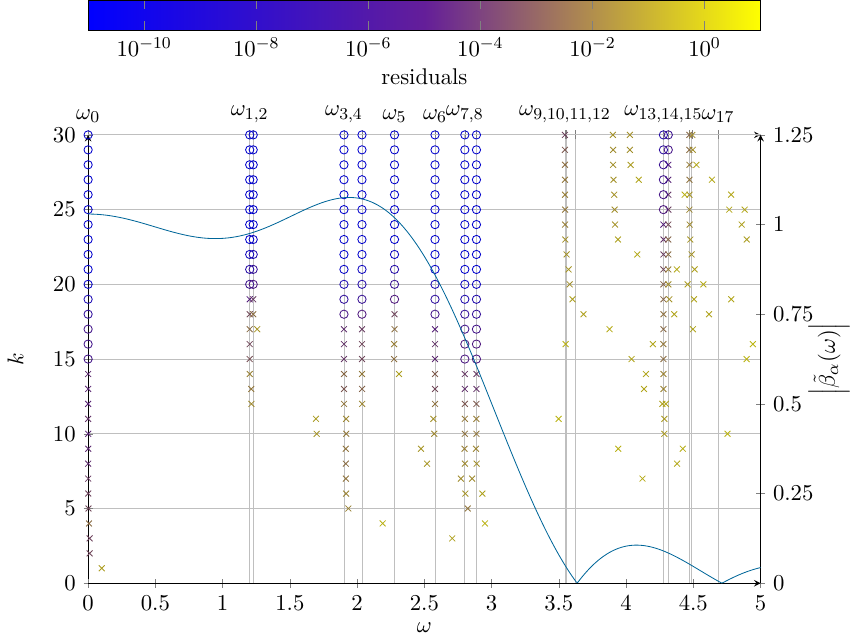}
  \caption{$\eigmi=0,\eigma=3, T=500\tau$}
  \label{fig:dumbbell_0_3_500}
  \end{subfigure}\\
  \begin{subfigure}{0.5\textwidth}
    \includegraphics[width=\textwidth]{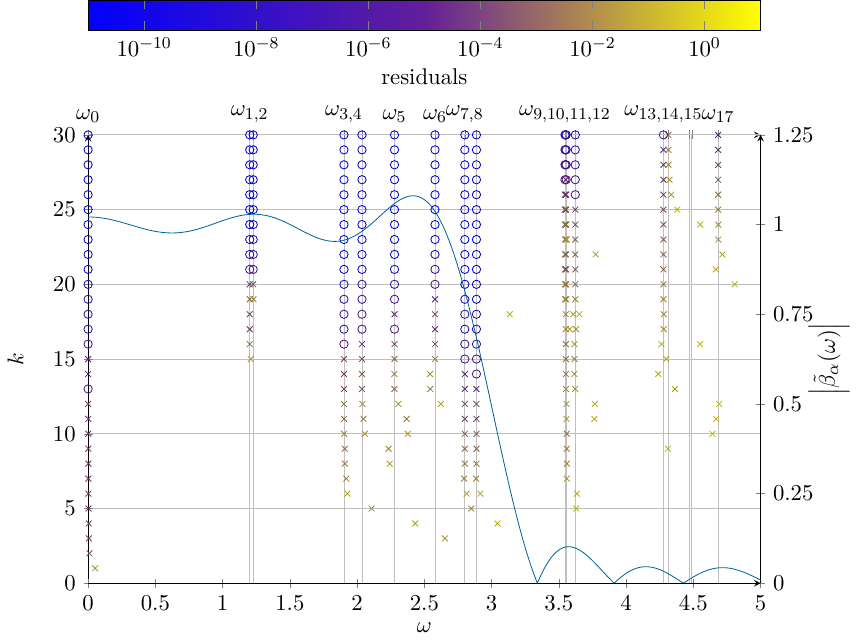}
  \caption{$\eigmi=0,\eigma=3, T=1000\tau$}
  \label{fig:dumbbell_0_3_1000}
  \end{subfigure}
  \begin{subfigure}{0.5\textwidth}
  \includegraphics[width=\textwidth]{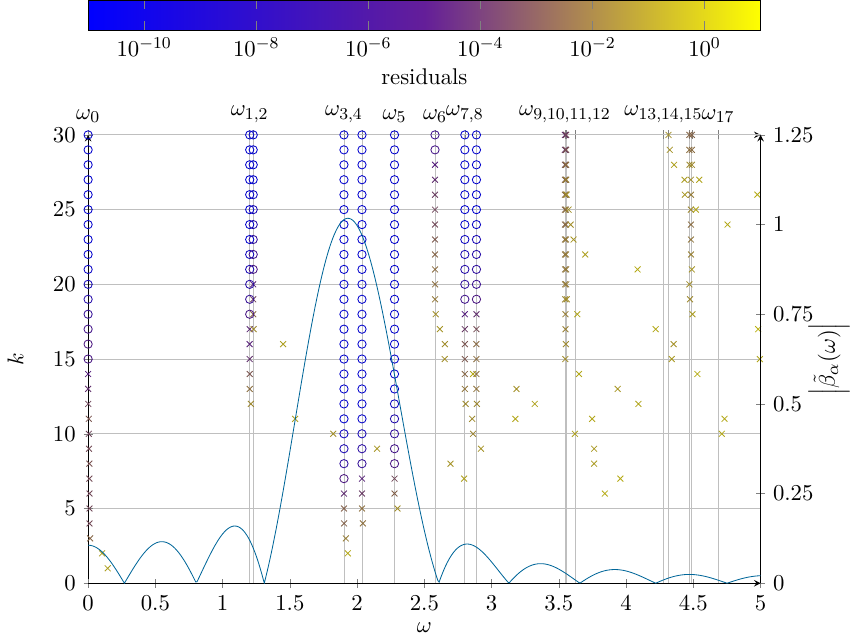}
  \caption{$\eigmi=1.6,\eigma=2.3, T=1000\tau$}
  \label{fig:dumbbell_1.6_2.3_1000}
  \end{subfigure}
  \caption{Convergence of the smallest eigenpairs for different weight function parameters $\eigmi,\eigma,T$ (see Figure \ref{fig:dumbbell_td_evs_qual} for a detailed explenation of the illustration).}
  \label{fig:dumbbell_lf}
\end{figure}

\begin{figure}
  \begin{subfigure}{0.5\textwidth}
    \includegraphics[width=\textwidth,trim={100 100 100 150},clip]{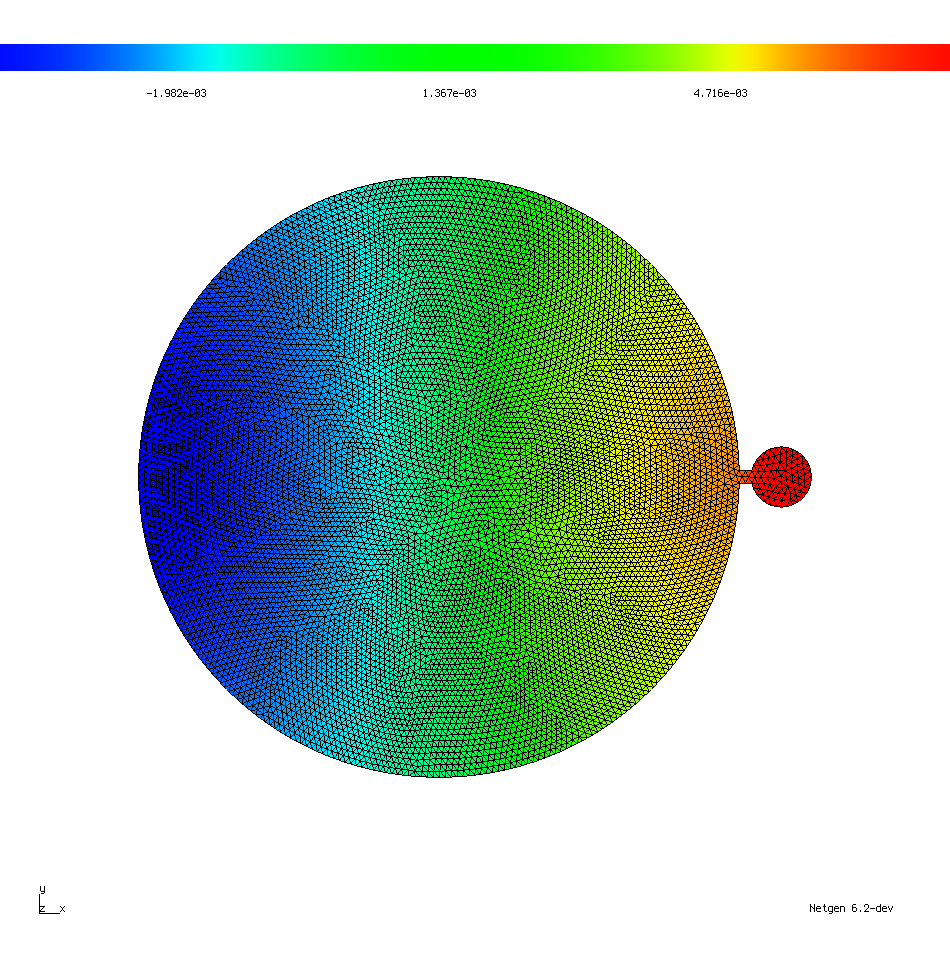}
    \caption{$\omega_{1}\approx1.2015$}
    \label{fig:2d_ef_1}
  \end{subfigure}
  \begin{subfigure}{0.5\textwidth}
  \includegraphics[width=\textwidth,trim={100 100 100 150},clip]{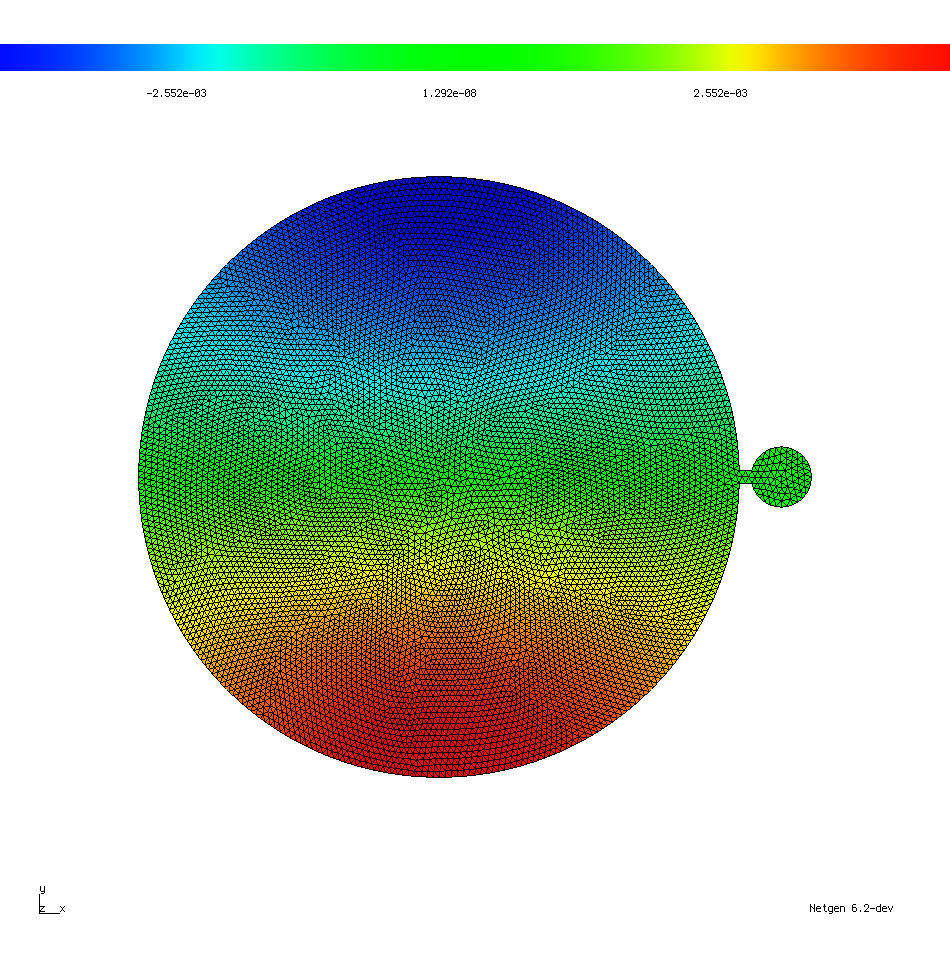}
    \caption{$\omega_{2}\approx1.2275$}
    \label{fig:2d_ef_2}
  \end{subfigure}
  \caption{Eigenfunctions corresponding to the resonances $\omega_1,\omega_2$ of the  2d example problem, where the green coloring is zero and red/blue marks positive/negative values.}
  \label{fig:2d_ef_low}
\end{figure}

We study the effect of choosing different parameters $\eigmi$,$\eigma$, and $T$ on the discrete filter and the approximation of eigenvalues. 
Figure \ref{fig:dumbbell_lf} shows the approximation of the same problem as before. In Figures \ref{fig:dumbbell_0_3_100}-\ref{fig:dumbbell_0_3_1000} the same interval $[\eigmi,\eigma]=[0,3]$ as before is combined with different numbers of time steps. In Fig.\ref{fig:dumbbell_1.6_2.3_1000} the interval is changed. Due to the few time steps $T=100\tau$, in Figure \ref{fig:dumbbell_0_3_100} the filter function does not have a good contrast for the wanted low frequencies and the unwanted higher frequencies. Thus convergence (except for $\omega=0$) is only achieved after $30$ Krylov steps. For more time steps (cf. Figures \ref{fig:dumbbell_0_3_500}, \ref{fig:dumbbell_0_3_1000} the contrast is better. This leads to convergence for eight resonances after $20$ Krylov steps. However the experiment also shows, that using more time steps does not necessarily lead to better results: The computational effort from Figure \ref{fig:dumbbell_0_3_500} to \ref{fig:dumbbell_0_3_1000} is doubled, due to the twice as many time steps. This would only pay off if the number of Krylov steps needed to compute the sought resonances  is more than halved. However using $1000$ time steps enables us to choose a narrow range of resonances (i.e., $\eigmi=1.6,\eigma=2.3$, in Figure \ref{fig:dumbbell_1.6_2.3_1000}). In this case the resonance $\omega_3$ is approximated well after 7 Krylov steps (7000 time steps in total) opposed to 18 Krylov steps (9000 time steps in total) in Figure \ref{fig:dumbbell_0_3_500}. 
To complete the  example Figure \ref{fig:2d_ef_low} shows the first two non-trivial eigenfunctions generated by the large ball. They are perturbations of the eigenfunctions the first non-trivial eigenvalue of multiplicity two of the closed ball. (cf. Section \ref{sec:dumbbell_desc}).

\begin{figure}
  \begin{subfigure}{0.5\textwidth}
    \includegraphics[width=\textwidth]{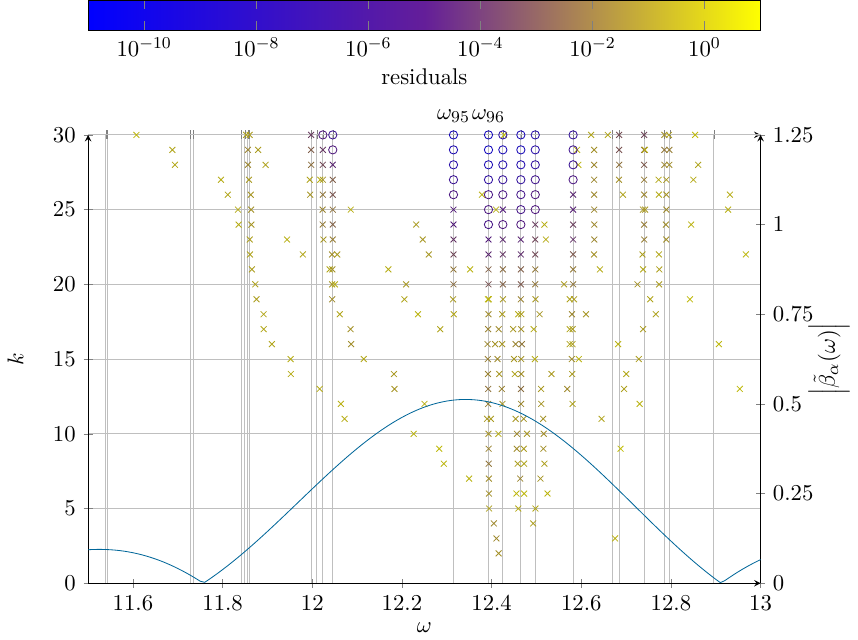}
  \caption{$T=1000\tau$}
  \label{fig:dumbbell_12.2_12.5_100}
  \end{subfigure}
  \begin{subfigure}{0.5\textwidth}
  \includegraphics[width=\textwidth]{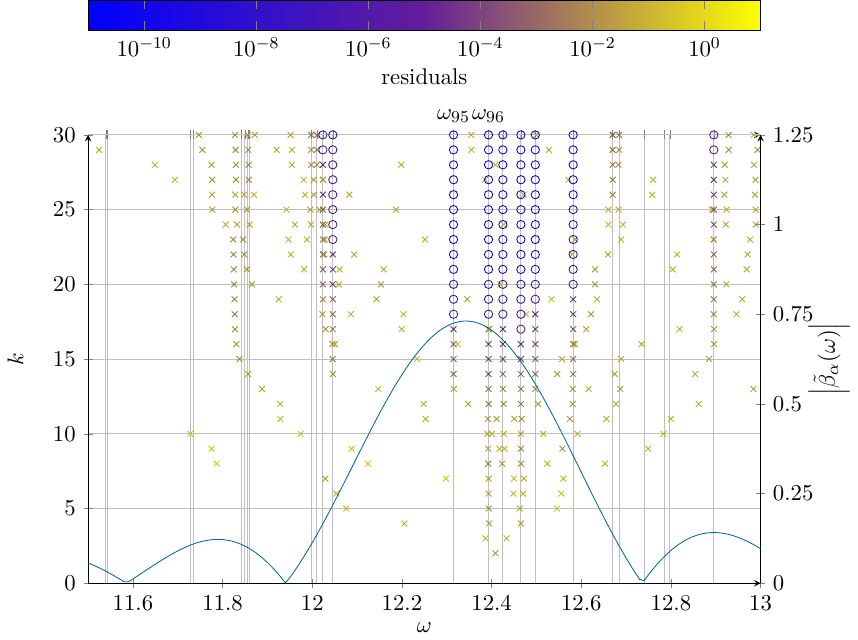}
  \caption{$T=1500\tau$}
  \label{fig:dumbbell_12.2_12.5_1500}
  \end{subfigure}
  \caption{Convergence of higher frequency eigenpairs for $\eigmi=12.2,\eigma=12.5$ and different end times $T$ (see Figure \ref{fig:dumbbell_td_evs_qual} for a detailed explanation of the illustration).}
  \label{fig:dumbbell_hf}
\end{figure}

To test the filtering of larger resonances 
we look for the resonance corresponding to the smallest non-trivial element of $\Lambda_r$ (cf. Figure \ref{fig:dumbbell_evs}). To this end we choose $\eigmi=12.2,\eigma=12.5$. The results in Figure \ref{fig:dumbbell_hf} show that the filtering of higher frequency resonances works as expected. However due to the fact that the spectrum of our problem gets denser in higher frequency ranges we have to choose a narrower peak of the filter and thus also more time steps.

\begin{figure}
  \begin{subfigure}{0.5\textwidth}
    \includegraphics[width=\textwidth,trim={100 100 100 150},clip]{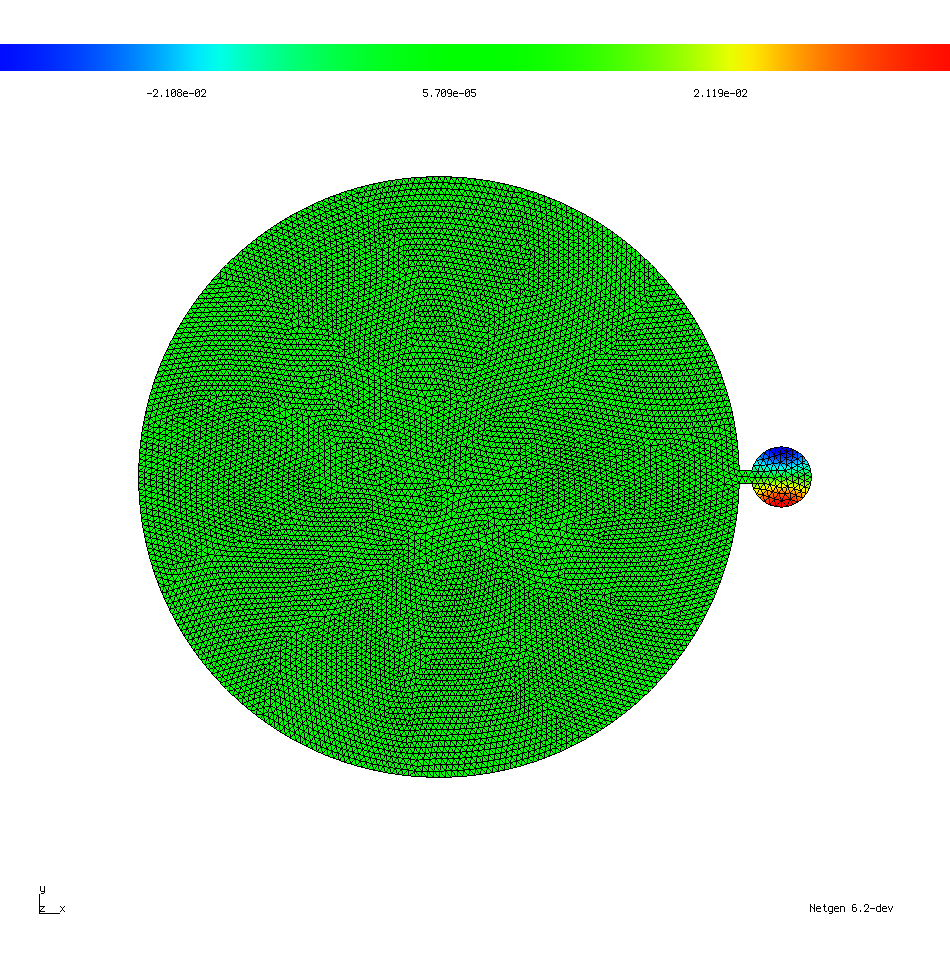}
    \caption{$\omega_{95}\approx12.3150$}
    \label{fig:2d_ef_95}
  \end{subfigure}
  \begin{subfigure}{0.5\textwidth}
  \includegraphics[width=\textwidth,trim={100 100 100 150},clip]{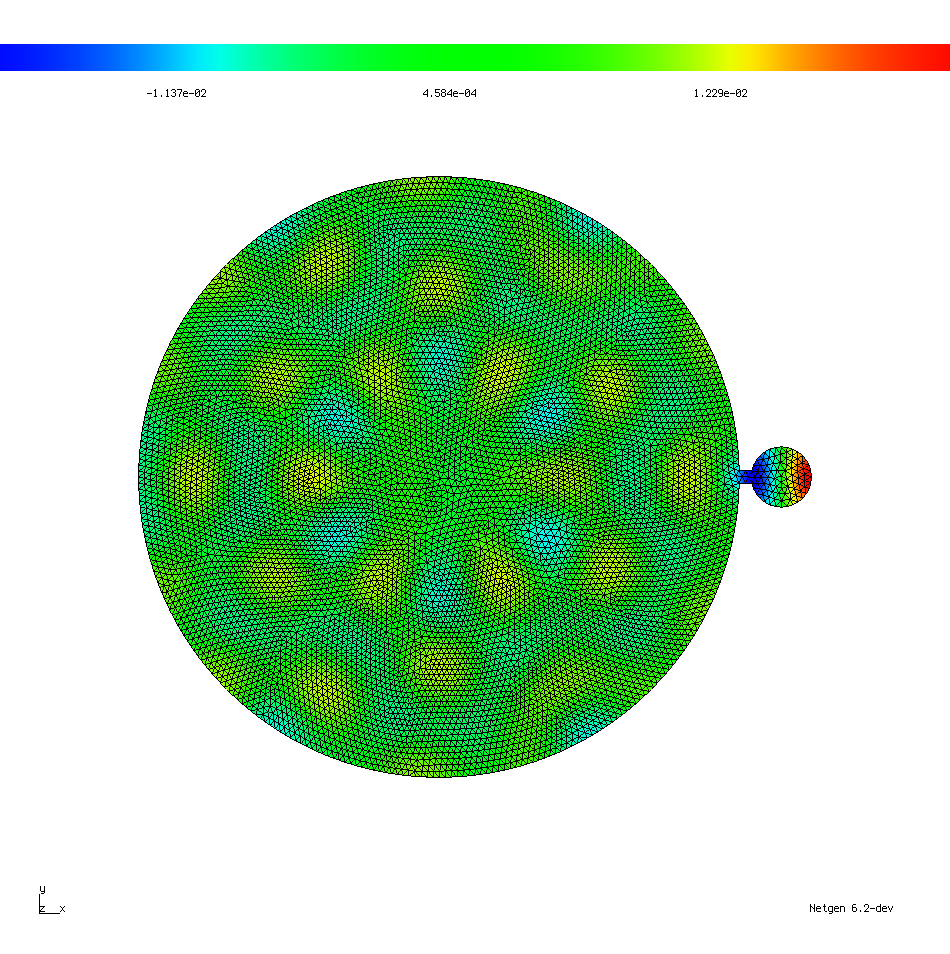}
    \caption{$\omega_{96}\approx12.3931$}
    \label{fig:2d_ef_96}
  \end{subfigure}
  \caption{Eigenfunctions corresponding to the resonances $\omega_{95},\omega_{96}$ of the  2d example problem.}
  \label{fig:2d_ef_high}
\end{figure}
 Figure \ref{fig:dumbbell_12.2_12.5_1500} also underlines the necessity to sort out the results based on the residuum, since even in the last Krylov step, there are non-converged resonances in between the converged ones. Finally, Figure \ref{fig:2d_ef_high} shows that the resonances around $12.3$ are in fact perturbations of the resonances of the closed small circle with multiplicity two. 

\FloatBarrier
\subsection{Large scale examples}
\label{sec:large_scale_numerics}
\begin{figure}[t]
  \centering
  \includegraphics[width=\textwidth, clip, trim={0 0 30 0}]{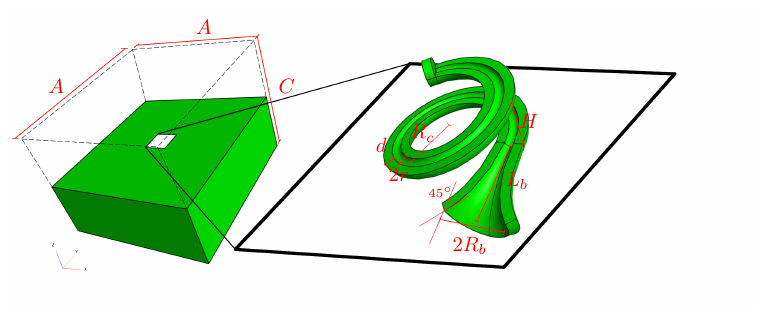}
  \caption{Cuts through the geometry of the horn and surroundings.}
  \label{fig:horn_geo}
\end{figure}
In this section we apply our method to large examples ($>10^6$ unknowns) where methods which require matrix inversion are no longer feasible on desktop computers, due to long factorization times and/or memory requirements. We demonstrate that our method is applicable on off-the-shelf computers.

\subsubsection{Description of the problem}
Our goal is to simulate a hunting horn in a closed room and to find the resonances corresponding to the lowest notes of the horn (i.e., the first few notes of the harmonic series).
The horn (cf. Figure \ref{fig:horn_geo}) consists of a coiled tube with inner radius $r$ and thickness $d$. The radius of the coil is $R_c$ and the winding number is $1.35$ and coil spacing $H$. The mouthpiece is modeled by Dirichlet boundary conditions on one end of the tube. The bell is modeled by the rotation of a spline along the axis with parameters $R_b,L_b$. 

The horn is enclosed in a cuboidal box with dimensions $A$, $B$,and $C$.
In all our experiments we fixed the parameters
\begin{align*}
  r=0.123,\qquad d=0.246,\qquad R_c=0.15, \qquad R_b=0.08,\qquad L_b=0.2,\qquad H=4r+2d
\end{align*}
of the horn and $A=B=3$, $C=2.4$ for the bounding box. To discretize the problem we use a mesh with mesh-size $h=d=0.0246$. We use a space of first order, mass lumped finite elements. This results in a problem of size $N\approx 1.245\cdot 10^6$. Again, a power iteration determines the time step $\tau \approx 0.000346$ to be stable.

Similar to the small scale experiments we expect the resonances of the resulting problem to be either perturbations of the resonances of the box given by
\begin{align*}
  \lambda_{i,j,k}&=\pi\sqrt{i^2/A+j^2/B+k^2/C},&i,j,k\in\setN_0,
\end{align*}
or perturbations of the resonances of the (closed) inner of the horn. We look for the latter, since they correspond to the musical notes which can be played on the horn.

  The example is chosen as a challenging task for eigenvalue solvers, since the sought resonance frequencies of the horn lie in a region where the background resonances of the box are already quite dense. Note however, that for the efficient computation of such resonances absorbing boundaries/layers like perfectly matched layers (PMLs) or infinite elements may be used.

\begin{figure}
  \begin{subfigure}{0.5\textwidth}
  \includegraphics[width=\textwidth]{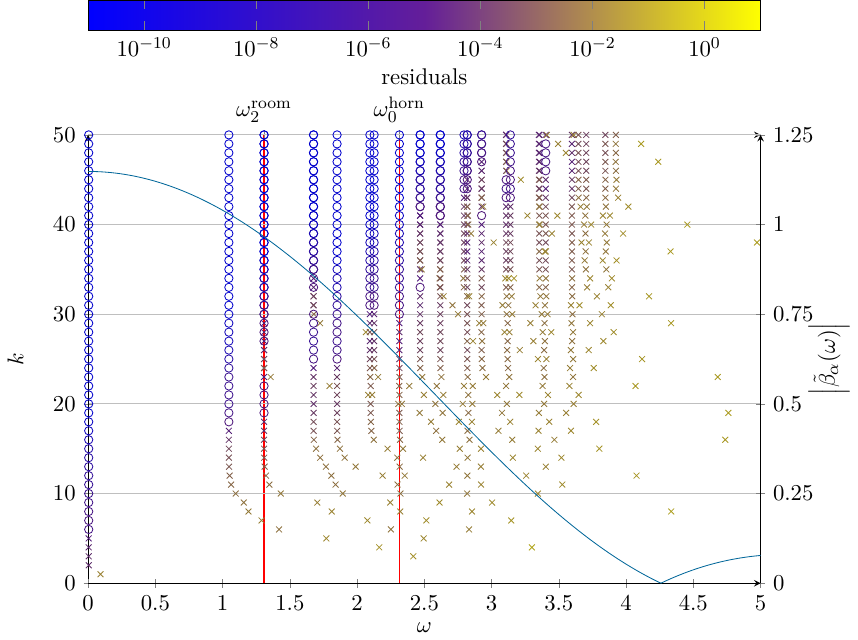}
    \caption{$T = 4000\tau, \eigmi=0,\eigma=2.7$}
    \label{fig:horn_low_fq}
  \end{subfigure}
  \begin{subfigure}{0.5\textwidth}
  \includegraphics[width=\textwidth]{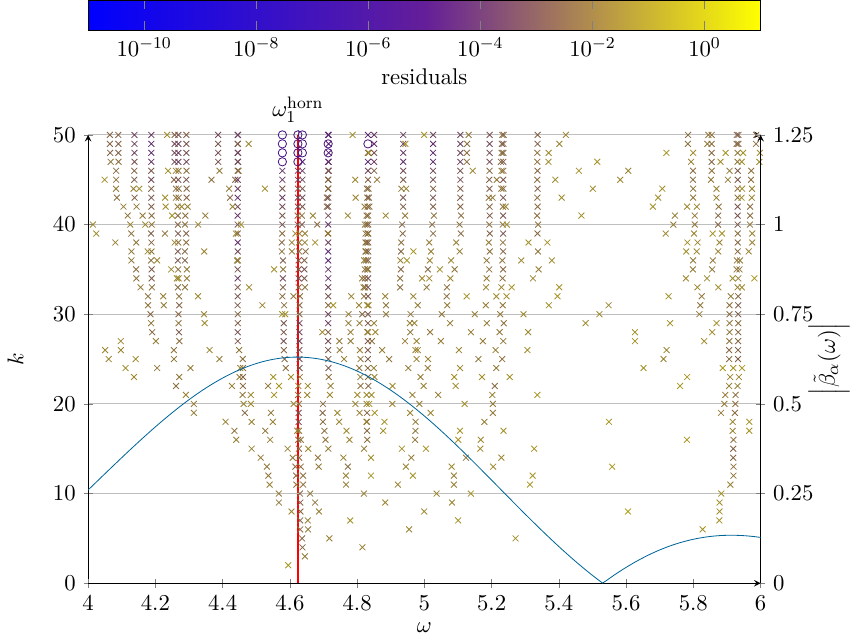}
    \caption{$T = 10000\tau, \eigmi=4.3,\eigma=4.9$}
    \label{fig:horn_mid_fq}
  \end{subfigure}
    \caption{Computed resonances of the horn domain for different filter parameters $T,\eigmi,\eigma$, where we look for the base harmonics of the horn $\omega^{\mathrm{horn}}_0,\omega^{\mathrm{horn}}_1$}
    \label{fig:horn}
\end{figure}

\begin{figure}
  \centering
  \begin{subfigure}{0.3\textwidth}
    \includegraphics[width=\textwidth, clip, trim = {250 30 250 100}]{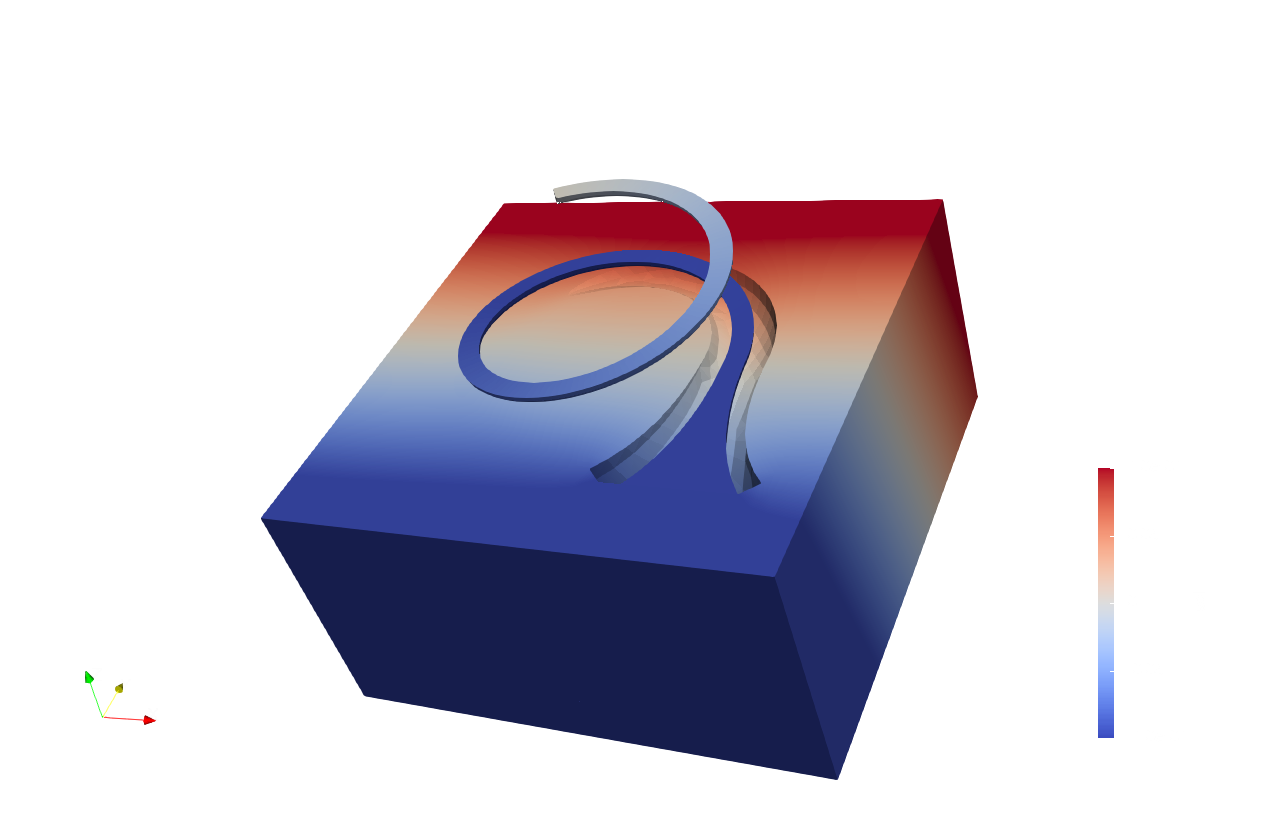}
    \caption{An exterior resonance $\omega^{\mathrm{room}}_2\approx 1.308$}
    \label{fig:horn_low_ef_outer}
  \end{subfigure}\hfill
  \begin{subfigure}{0.3\textwidth}
  \includegraphics[width=\textwidth, clip, trim = {250 30 250 100}]{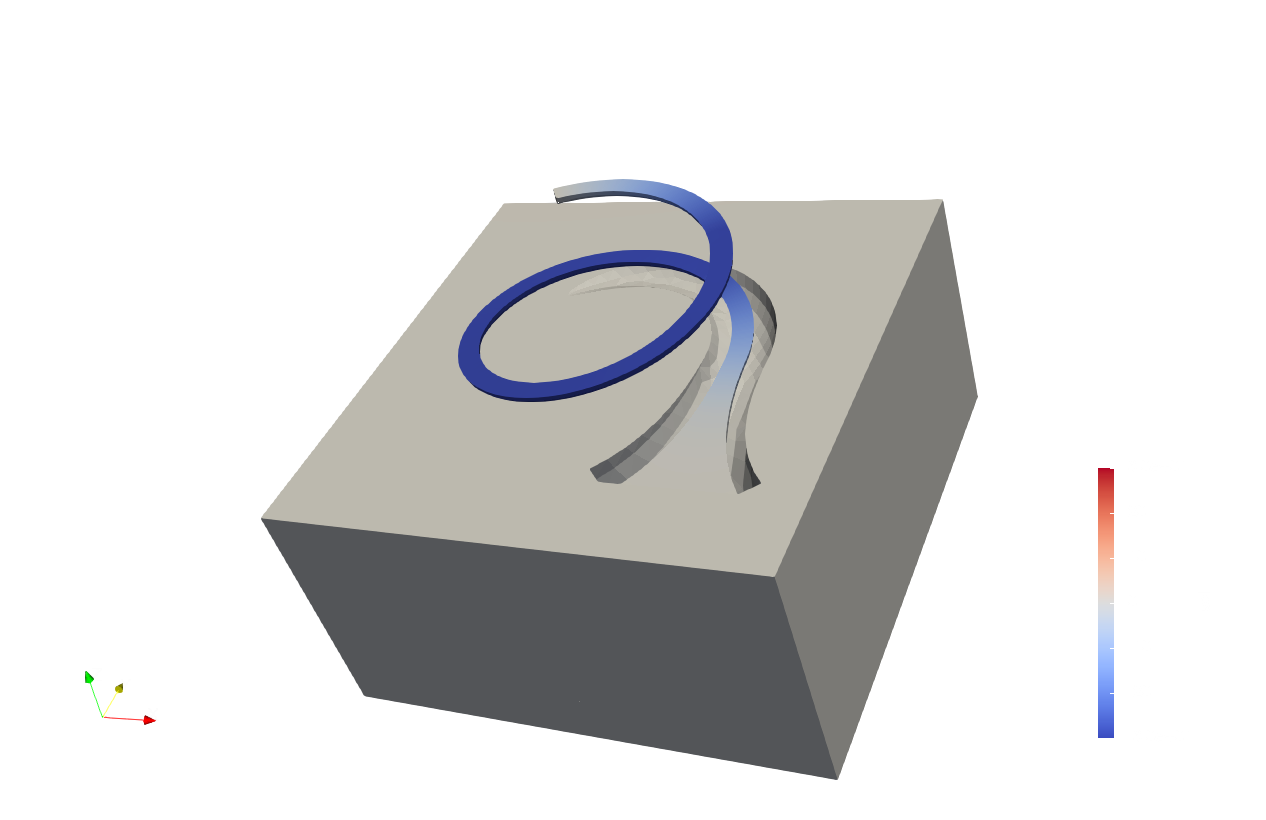}
    \caption{Base resonance at $\omega^{\mathrm{horn}}_0 \approx 2.315$}
    \label{fig:horn_low_ef_inner}
  \end{subfigure}\hfill
  \begin{subfigure}{0.3\textwidth}
  \includegraphics[width=\textwidth, clip, trim = {250 30 250 100}]{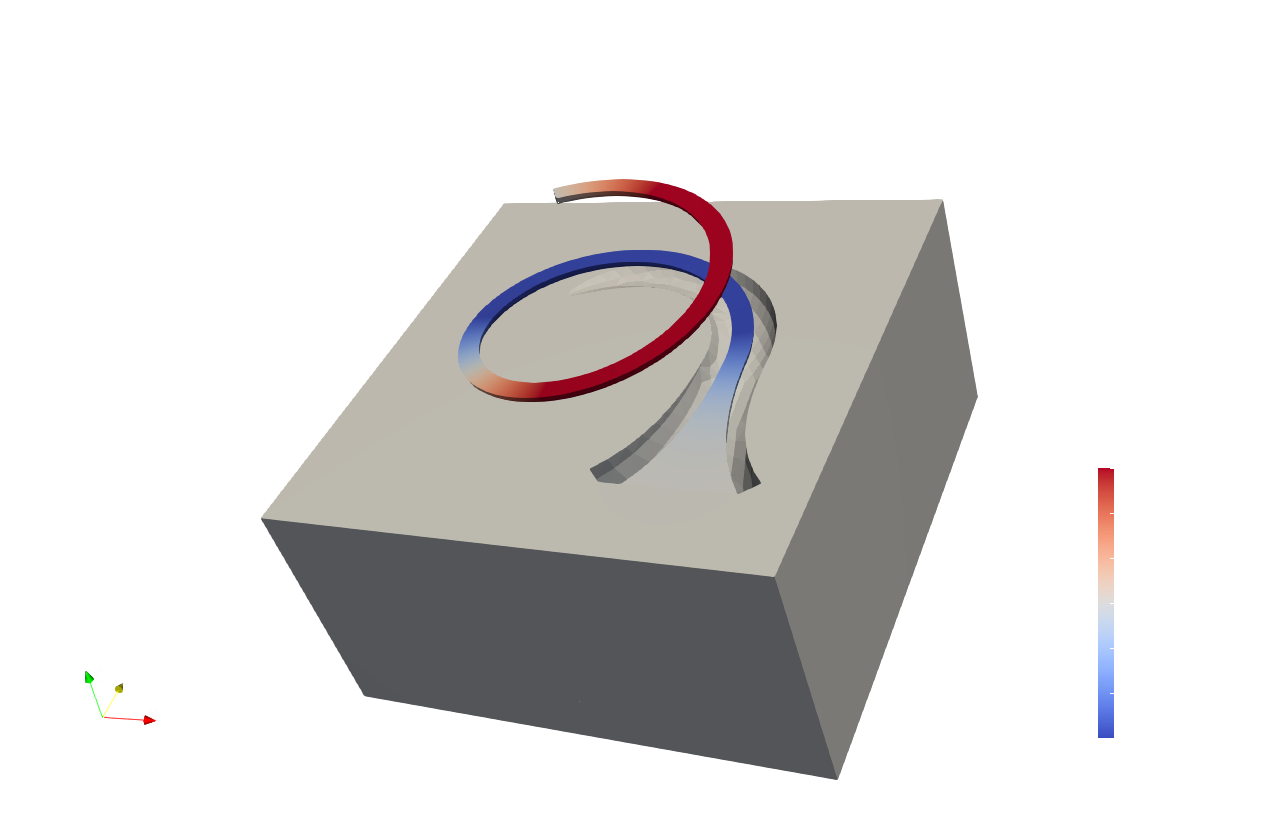}
    \caption{First harmonic at $\omega^{\mathrm{horn}}_1 \omega\approx 4.624$}
    \label{fig:horn_mid_ef}
  \end{subfigure}
    \caption{Eigenfunctions of the horn problem corresponding to the frequencies $\omega^{\mathrm{room}}_2,\omega^{\mathrm{horn}}_0,\omega^{\mathrm{horn}}_1$ from Figure \ref{fig:horn}. The colors blue and red correspond to higher/lower values while grey corresponds to zero. Note that not the whole domain, but merely a cutout (cf., Figure \ref{fig:horn_geo}) is shown.}
    \label{fig:horn_ef}
\end{figure}

As before we use different sets of parameters for $T,\eigmi,\eigma$ to search in different frequency ranges for eigenfrequencies. Using $T=4000\tau,\eigmi=0,\eigma=2.7$ we find the base mode of the horn at $\omega^{\mathrm{horn}}_0\approx 2.315$ (cf. Figures \ref{fig:horn_low_fq}, \ref{fig:horn_low_ef_inner}). Thus we look for the second harmonic at approximately $2\omega^{\mathrm{horn}}_0$. Figure \ref{fig:horn_mid_fq} shows the result for $T=10000,\eigmi=4.3,\eigma=4.9$. The spectrum is already very dense in this region, still our algorithm manages to converge for a few resonances. Indeed we find the second harmonic of our horn at $\omega^{\mathrm{horn}}_1\approx4.624$ (cf. Figure \ref{fig:horn_mid_ef}).

\begin{rmk}
  If we assume the dimensions of our horn to be given in meters and a speed of sound of $343m/s$ we obtain that the first two notes which can be played have frequencies of 
  \begin{align*}
    343\frac{\omega^{\mathrm{horn}}_0}{2\pi}&\approx 126.38 Hz,&
    343\frac{\omega^{\mathrm{horn}}_1}{2\pi}&\approx 252.43 Hz.&
  \end{align*}
  In musical notation this corresponds to notes a little higher than a great and small B respectively. This is a reasonable result given the fact that our horn has a total length of $1.479m$, which is in a similar range as a trumpet tuned in B flat (usually a length of approximately $1.485m$).
\end{rmk}

\FloatBarrier

\section{Generalizations and extensions}
\label{sec:generalizations}
In Section \ref{sec:basic_example} and the corresponding numerical experiments in Section \ref{sec:numerics} we focused on applying our ideas to an example composed of the following components:
\begin{description}
  \item{\bf discrete problem: }{The generalized matrix eigenvalue problem \eqref{def:basisMatrixEVP} stems from a finite element discretization of 
  a Laplacian eigenvalue problem using mass lumping.}
  \item{\bf time-stepping: }{The discrete filter function defined in \eqref{eq:fullydiscfilter} is based on the Verlet time-stepping of the time-domain problem \eqref{eq:waveeq} corresponding to \eqref{def:basisMatrixEVP} and a suitable} 
  \item{\bf weight function: }{The weight function is based on an approximation of the inverse Fourier transform of an indicator function (cf. Equations \eqref{eq:filterf} and \eqref{eq:weightf}).}
  \item{\bf eigenvalue solver: }{To approximate the eigenvectors of the auxiliary eigenvalue problem \eqref{def:ArnoldiEWP} we construct a basis of the Krylov space of the operator $C$ by iterative application and orthonormalization.}
\end{description}
In fact, all of these components can replaced by suitable alternatives. In the following we discuss such extensions of the method, as well as the limits of the current approach. 

\subsection{Discrete Problem}
\label{sec:gen_disc}

The main requirement for the method to work is the fact that the matrix eigenvalue problem \eqref{def:basisMatrixEVP} corresponds to a stable time-domain problem of the form \eqref{eq:waveeq}. This is certainly the case if the matrices $M$ and $S$ 
are symmetric and positive definite and positive semi-definite respectively. This is fulfilled for suitable discretizations of wave-type equations where $S$ is the discrete representation of any elliptic, second order differential operator. Examples include the $\curl\curl$-operator (for linear Maxwell equations) or the elastic operator.

The main requirement for the method to be efficient is the fact that the according time-domain problem can be efficiently approximated. For explicit time-stepping methods (cf. Section \ref{sec:gen_time_stepping}) this is the case if the inverse of the mass matrix $M$ can be applied efficiently. Apart from mass lumping techniques this is also the case for finite difference methods, discontinuous Galerkin approaches (see e.g., \cite{hesthavenDG}) or cell methods \cite{cell_method_lumped,cell_method}.

Going further one could even drop the assumption that the problem is self-adjoint (Hermitian) if there is a stable time-domain counterpart. One example could be resonance problems in open systems, which are discretized using perfectly matched layers (cf. \cite{Berenger:94}). However,  the choice of the weight function $\wei$ is not a-priori clear if complex resonances exist.


\subsection{Time-stepping}
\label{sec:gen_time_stepping}
Following Remark \ref{rmk:timestep} the Verlet time-stepping  \eqref{eq:2stepmethod} can be easily replaced  by any other explicit time-stepping method, if the eigenvectors of \eqref{def:basisMatrixEVP} are still eigenvectors of \eqref{def:ArnoldiEWP}. For first order formulations one could e.g., use leap frog time-stepping or high-order variants thereof to construct a similar method. 
Time-stepping schemes where the correspondence of eigenvectors as described above is not clear any more include local time-stepping \cite{Grote:09} or locally implicit methods \cite{Piperno:06}. These methods have been proven to significantly reduce the CFL condition on the time-step size which is a major limiting factor in computational efficiency. However, they do not fit directly into the concept presented in this paper.

\subsection{Weight function}
\label{sec:gen_weight}
The weight function chosen in \eqref{eq:weightf} is based on the inverse Fourier transform of a characteristic function. As an alternative to the characteristic function one could choose more localized and/or smoother functions (e.g., a triangle impulse or a Gaussian peak) to construct the weight function $\wei$. 
%
In particular, the weight function $\wei$ could be adapted such that in regions, where an essential spectrum or an eigenvalue with very high multiplicity is expected, the absolute values of the corresponding discrete filter $\tilde \beta_\wei$ is small. 
%
%
An example where this would be useful is e.g., the infinite dimensional kernel of the $\curl\curl$-operator.
\subsection{Eigenvalue solver}
To resolve issues with multiple eigenvalues a block version of the Krylov space approach can be used (i.e., starting with an orthonormalized set of random vectors and also orthonormalizing the resulting vectors in each step). Alternatively instead of considering a growing Krylov space in each step one could fix the dimension which leads to a method corresponding to a FEAST \cite{FEAST} algorithm with a time-domain filter.

\section{Conclusion}
\label{sec:conclusion}

We have presented an eigenvalue solver for large scale eigenvalue problems originating from finite element discretizations of Laplacian eigenvalue problems. Existing explicit time-domain solvers are combined with a weight function in order to focus a Krylov method to eigenvalues, which might be non-extremal and/or clustered. It is straightforward to show, that the method converges to those eigenvalues, which are mapped by the discrete filter $\tilde \beta_\wei$ defined in \eqref{eq:fullydiscfilter} to eigenvalues with largest absolute value of the auxiliary eigenvalue problem. A complete error analysis would need to control the absolute values of the filtered sought and unsought eigenvalues, which lies out of the scope of this paper. 

Note, that the discrete filter depends on the weight function $\wei$, the number and size of time-steps, and 
on the time-stepping scheme. The filter function can be computed in a preprocessing step with negligible computational costs. Hence, an optimization of the method parameters for a given problem is easy to carry out experimentally.

\bibliographystyle{alpha} 
\bibliography{bibliography}
\end{document}